\documentclass{article}
\usepackage[utf8]{inputenc}
\usepackage[russian, english]{babel}
\usepackage{hyperref}
\usepackage{graphicx}
\usepackage[T1,T2A]{fontenc}
\usepackage{xcolor}
\usepackage{gensymb}
\usepackage{amsmath,amssymb,amsthm,amsfonts}
\newtheorem{theorem}{Theorem}

\newcommand{\gf}{\mathfrak{g}}

\begin{document}

\title{Probability measure near the boundary of tensor power decomposition for $\mathfrak{so}_{2n+1}$}

\author{Anton Nazarov$^{1,2}$, Viktoria Chizhikova$^{1,3}$\\
{\small$^{1}$Department of Physics, St. Petersburg State University,} \\
{\small  Ulyanovskaya 1, 198504 St.~Petersburg, Russia}\\
{\small$^{2}$email:antonnaz@gmail.com}\\
{\small$^{3}$email:vika\_chizh@mail.ru}}
\maketitle

\begin{abstract}
    Character measure is a probability measure on irreducible representations of a semisimple Lie algebra. It appears from the decomposition into irreducibles of tensor power of a fundamental representation. In this paper we calculate the asymptotics of character measure on representations of $\mathfrak{so}_{2n+1}$ in the regime  near the boundary of weight diagram. We find out that it converges to a Poisson-type distribution.
\end{abstract}

\section{Introduction}

Consider a simple finite-dimensional complex Lie algebra $\gf$ of rank $n$ and its irreducible finite-dimensional highest-weight representation $L^{\lambda}$. Denote
simple roots of $\gf$ by $\alpha_{1},\dots, \alpha_{n}$ and fundamental weights by $\omega_{1},\dots \omega_{n}$:
$(\alpha_{i},\omega_{j})=\delta_{ij}$. 
Let us take a fundamental representation $L^{\omega}$ of $\gf$, $\omega\in\{\omega_1, ... , \omega_n\}$,
and consider its tensor power. Tensor power of $L^{\omega}$ is a completely reducible representation and can be decomposed as:
\begin{equation}
  \left(L^{\omega}\right)^{\otimes N}=\bigoplus_\lambda M_\lambda^{\omega,N} L^\lambda,
\end{equation}
where $M_\lambda^{\omega,N}$ is a multiplicity of $L^\lambda$.
The sum is taken over all irreducible components of the tensor product. This means that for the dimensions holds the equality
\begin{equation}
  \dim\left(\left(L^{\omega}\right)^{\otimes N}\right)=\sum_\lambda M_\lambda^{\omega,N} \dim( L^\lambda).
\end{equation}

This formula can be used to introduce Plancherel-type probability measure on the set of
dominant integral weights: 
\begin{equation}
  \label{eq:3}
  P_\lambda^N=\frac{M_\lambda^{\omega,N} \dim L^\lambda}{(\dim L^{\omega})^N}.
\end{equation}

We can generalize this formula by replacing $\dim(L^\lambda)$, which is equal to $\chi_{L^\lambda}(0)$, by a character $\chi_{\lambda}(e^{t})$ of $L^\lambda$ taken at any item $e^t$ from maximal torus of the Lie group or item $t$ from the Cartan subalgebra $\mathfrak{h}$ of the Lie algebra $\mathfrak{g}$. The character of a representation $L^\lambda$ of a Lie algebra is equal to 
\begin{equation}
  \label{eq:4}
 \chi_\lambda(e^t)=\sum_{\mu\in \mathcal{N}(\lambda)} \dim V_\mu \cdot e^{(\mu,t)}.
 \end{equation}
Representation $L^\lambda=\bigoplus_\mu V_\mu$ is decomposed into direct sum of weight subspaces and $\dim V_\mu=m^\lambda_\mu$ is the multiplicity of weight $\mu$ in $L^\lambda$. $\mathcal{N}(\lambda)$ is the weight diagram.
So we get the character measure:
\begin{equation}
  \label{eq:5}
P_\lambda^N(t)=\frac{M_\lambda^{\omega,N} \chi_\lambda(e^t)}{(\chi_{\omega}(e^t))^N}.
\end{equation}

 The measures \eqref{eq:3},\eqref{eq:5} in the limit $n\rightarrow\infty$ show the statistics of representations of infinite-dimensional algebras; in the limit $N\rightarrow\infty$ -- the statistics of infinite tensor products. Plancherel measure for permutation group $S_{n}$ is also connected to Ulam’s problem on the length of the maximal increasing subsequence in a uniform random sequence. Its central limit regime is well studied and is known as limit shape of Young diagrams \cite{vershik}. The same result is obtained for measure \eqref{eq:3} if $\mathfrak{g}=\mathfrak{sl}_{n}$ in the limit $N,n\to \infty$ in the paper \cite{kerov1986}, where the case $N\to\infty,n=\mathrm{const}$ was also considered. The case $\mathfrak{g}=\mathfrak{so}_{2n+1}, N\to\infty, n=\mathrm{const}$ was studied in \cite{nazarov2018limit} and \cite{tate2003lattice}. The central limit regime and the regime of large deviations for character measure \eqref{eq:5}, with $n$ fixed, $N\rightarrow\infty$, were considered in \cite{postnova2020multiplicities} by O.Postnova and N.Reshetikhin. They also suggested that in the regime near the boundary the character measure will converge to some Poisson type process. They demonstrated that for $\mathfrak{g}=\mathfrak{sl}_2$ it converges to Poisson distribution; in this paper we generalize this result to all $\mathfrak{so}_{2n+1}$ algebras.

\begin{theorem}
Consider Lie algebra $\mathfrak{so}_{2n+1}$  with root system  $B_n$, $n$ is fixed.
Let us take the tensor power $N$ of its irreducible representation with the highest weight $\lambda=\omega_n$ (it is known as last fundamental representation or spinor representation). Consider the character measure  $P_\lambda^N(t)$ given by the formula \eqref{eq:5}, taken at $e^{t}$, where $t$ is an element of Cartan subalgebra of $\mathfrak{so}_{2n+1}$.

Then as  $N\rightarrow \infty$ and $t\rightarrow\infty$ such that $\Theta_i=N e^{-2t_i}$ is fixed, the probability measure $P_\lambda^N(t)$ converges pointwise to
\begin{equation}
  \label{eq:1}
  P_s(\Theta)=\prod_{i<j}^n(-s_i-i+s_j+j) \cdot  \prod_{k=1}^{n} \frac{  \Theta_k^{\frac{s_1+...+s_n}{n}} \cdot e^{-\Theta_k}}{(s_{k}+k-1)!}\chi^{\mathfrak{sl}_{n}}_\gamma(e^\tau),
\end{equation}
where $N-2s_i$ are coordinates of the highest weight $\lambda$, $\chi^{\mathfrak{sl}_{n}}_\gamma(e^{\tau})$ is a character of the  $\mathfrak{sl}_{n}$-subalgebra representation of the highest weight $\gamma$ with the coordinates $\gamma_i=\frac{2s_1+...+2s_n}{n}-2s_i$ and $\tau$ is an element of Cartan subalgebra in $\mathfrak{sl}_{n}\subset \mathfrak{so}_{2n+1}$ with the coordinates $\tau_i=-t_i+t_{i+1}$. 
\end{theorem}

In the next section we present the proof of this theorem and discuss its meaning and possible generalizations in Conclusion. 
\section{Proof of the theorem}

\subsection{Coordinates and multiplicities}
We will use the standard orthogonal basis such that all the elements of Cartan subalgebra would be diagonal. In this basis simple roots are
\begin{equation}
  \label{eq:6}
\alpha_i=e_i-e_{i+1}, \; i=1, .... , n-1; \;\;\; \alpha_n=e_n,
\end{equation}
and fundamental weights are
\begin{equation}
  \label{eq:7}
\omega_1=e_1; \;\; \omega_2=e_1+e_2; \;\; ...\;\; \omega_{n-1}=e_1+...+e_{n-1};  \;\; \omega_n=\frac{1}{2}(e_1+...+e_n) .
\end{equation}
For the highest weight $\lambda$ we denote its coordinates by $\widetilde{\lambda_i}$ so that $\lambda=\sum_{i=1}^{n}\widetilde{\lambda_i}e_i$. For convenience we rescale the coordinates: $\lambda_i=2\widetilde{\lambda_i}$. The multiplicity formula \cite{kulish2012multiplicityBn} can be written easily using the coordinates shifted by the Weil vector $\rho=\omega_1+...+\omega_n$, $a_i=\lambda_i+\rho_i$, where $\rho_i=2(n-i)+1$ are the coordinates of $\rho$:
\begin{equation}
  \label{eq:8}
M^{\omega_n,N}_\lambda=\prod_{k=0}^{n-1} \dfrac{(N+2k)!}{2^{2k} \left(\frac{N+a_{k+1}+2n-1}{2}\right)!\left(\frac{N-a_{k+1}+2n-1}{2}\right)!} \prod_{l=1}^na_l \prod_{i<j}^n(a_i^2-a_j^2)
\end{equation}
We study regime near the boundary, which means $\lambda$ should be not far from the highest weight of $(L^{\omega_n})^{\otimes N}$ representation. So we take
\begin{equation}
  \label{eq:9}
\lambda_i=N-2s_i, \;\;\;  s_i=\mathcal{O} (1).
\end{equation}
Let us rewrite the formula  \eqref{eq:8} using the parameters $\{s_i\}$:
\begin{equation}
  \label{eq:10}
M^{\omega_n,N}_\lambda=\prod_{k=0}^{n-1} \dfrac{(N+2k)!}{2^{2k} \left(N+\frac{-2s_{k+1}+\rho_{k+1}+2n-1}{2}\right)!\left(\frac{2s_{k+1}-\rho_{k+1}+2n-1}{2}\right)!} \cdot\widetilde{M},
\end{equation}
where we have denoted by $\tilde{M}$ the following expression:
\begin{equation}
  \label{eq:11}
\widetilde{M}=\prod_{l=1}^n (N-2s_l+\rho_l) \prod_{i<j}^n(-2s_i+\rho_i+2s_j-\rho_j)(2N-2s_i+\rho_i-2s_j+\rho_j).
\end{equation}
Substituting the coordinates $\{\rho_{i}\}$ and using the identity $\prod_{k=0}^{n-1}2^{2k}=2^{n(n-1)},$ we get
\begin{equation}
  \label{eq:12}
M^{\omega_n,N}_\lambda=\frac{1}{2^{n(n-1)}}\prod_{k=1}^{n} \dfrac{(N+2(k-1))!}{ (N-s_k+2n-k)!(s_k+k-1)!} \cdot\widetilde{M},
\end{equation}
with
\begin{equation}
  \label{eq:13}
\widetilde{M}=\prod_{l=1}^n (N-2s_l+2n-2l+1)\cdot\prod_{i<j}^n(-2s_i-2i+2s_j+2j)(2N-2s_i-2i-2s_j-2j+4n+2).
\end{equation}

First we write down the leading asymptotic in $N$ for $\widetilde{M}$:
\begin{equation}
  \label{eq:14}
\widetilde{M}=N^n\cdot (2N)^{\frac{n(n-1)}{2}} \cdot\prod_{i<j}^n(-2s_i-2i+2s_j+2j)\left(1+\mathcal{O} \left(  \frac{1}{N} \right) \right).
\end{equation}

Then we derive the leading asymptotic in $N$ for another factor in the expression \eqref{eq:12}:
$$\dfrac{(N+2(k-1))!}{ (N-s_k+2n-k)!}=N^{N+2k-2-N+s_k-2n+k}\left(1+\mathcal{O} \left(  \frac{1}{N} \right) \right)=N^{s_k+3k-2n-2} \left(1+\mathcal{O} \left(  \frac{1}{N} \right) \right).$$

The whole expression is:
\begin{multline}
  \label{eq:2}
  M^{\omega_n,N}_\lambda=\frac{1}{2^{n(n-1)}}\prod_{k=1}^{n} \left(\dfrac{N^{s_k+3k-2n-2}}{ (s_k+k-1)!} \right)\cdot N^n\cdot (2N)^{\frac{n(n-1)}{2}} \cdot\prod_{i<j}^n(-2s_i-2i+2s_j+2j)\left(1+\mathcal{O} \left(  \frac{1}{N} \right) \right)=\\
  =N^{-\frac{3}{2}n(n+1)} \prod_{k=1}^{n} \frac{  N^{s_k+3k}}{(s_{k}+k-1)!} \cdot\prod_{i<j}^n(-s_i-i+s_j+j)\left(1+\mathcal{O} \left(  \frac{1}{N} \right) \right).
\end{multline}

Substitute $\prod_{k=1}^n N^{3k}=N^{\frac{3}{2}n(n+1)}$ and obtain the asymptotic for the multiplicities:
\begin{equation}
  \label{eq:15}
M^{\omega_n,N}_\lambda= \left[\prod_{k=1}^{n} \frac{N^{s_k}}{(s_{k}+k-1)!} \cdot\prod_{i<j}^n(-s_i-i+s_j+j)\right]\cdot\left(1+\mathcal{O} (  \frac{1}{N} ) \right).
\end{equation}
Now we need to derive the asymptotic of characters. 

\subsection{Asymptotics of characters}
\label{sec:asympt-char}
First we consider the character of $L^{\omega_n}$. Weights of $L^{\omega_{n}}$ are located at the vertices of an $n$-dimensional cube. Their coordinates are all combinations of 1 and -1, thus for the character in the $N$-th power we have
\begin{equation}
  \label{eq:16}
(\chi_{\omega_n}(e^t))^N=(e^{t_1+...+t_n}+e^{-t_1+t_2+...+t_n}+e^{t_1-t_2+...+t_n}+...+e^{-t_1-...-t_n})^N.
\end{equation}
Rewriting this expression using the parameters $\{\Theta_{i}\}$, in the limit we can neglect the terms like $e^{-t_i}$ as $e^{-t_i}=\mathcal{O}\left(\frac{1}{N}\right)$:
\begin{multline}
  \label{eq:20}
\chi_{\omega_n}(e^t)= e^{t_1+...+t_n}+e^{-t_1+t_2+...+t_n}+...+e^{t_1+...-t_i+...+t_n}+...+e^{t_1+...+t_{n-1}...-t_n}+\\ +e^{t_1+...+t_n}\cdot \mathcal{O} \left(  \frac{1}{N^2} \right)=e^{t_1+...+t_n}\left(1+e^{-2t_1}+...+e^{-2t_n}+O \left(  \frac{1}{N^2} \right)\right).
\end{multline}

So, using the new parameters we get:
\begin{equation}
  \label{eq:17}
  (\chi_{\omega_n}(\Theta))^N=\left(e^{t_1+...+t_n}\left(1+\frac{\Theta_1+...+\Theta_n}{N}+\mathcal{O} \left(  \frac{1}{N^2} \right)\right) \right) ^N,
\end{equation}
which is equal to 
\begin{equation}
  \label{eq:18}
(\chi_{\omega_n}(\Theta))^N= e^{(t_1+...+ t_n)N}\cdot e^{\left(\Theta_1+...+\Theta_n\right)\left(1+\mathcal{O} \left(  \frac{1}{N} \right)\right) }.
\end{equation}
 
Then we calculate the leading contribution to the character of $L^{N-2s}$. Remember that terms of the same order as $e^{-t_i}$ could be neglected. What weights contribute to this asymptotic?
These are the weights with the maximal sum of coordinates, denote this set of weights by $\Omega\subset \mathcal{N}(\lambda)$. These weights belong to the plane that contains the highest weight $\lambda$ and is orthogonal to $(1, ... , 1)$. All these weights are obtained from the highest weight by subtracting the simple roots that lie in this plane.
Simple roots lying in this plane are all the roots except $\alpha_n$, since for $k\in\{1,\dots,n-1\}$ we have:
\begin{equation}
  \label{eq:19} 
(1, 1, ... , 1)\cdot\alpha_k=\frac{1}{2}(1, ... , 1)\cdot(0, ..., 1, -1, ..., 0)=0.
\end{equation}
These roots make up the root system $A_{n-1}$, which corresponds to a regular subalgebra $\mathfrak{sl}_n$ in $\mathfrak{so}_{2n+1}$. Since the Weyl group $W_{\mathfrak{sl}_n}$ of the subalgebra  $\mathfrak{sl}_n$ is a subgroup of the Weyl group $W_{ \mathfrak{so}_{2n+1}}$ and leaves the plane invariant, the set $\Omega$ is also $W_{\mathfrak{sl}_n}$-invariant.  The weight subspaces for the weights in $\Omega$ are obtained from the highest weight vector $v_{\lambda}$ by a repeated action of the lowering generators $\{E^{-\alpha_{k}}, k=1,\dots,n-1\}$ corresponding to the $\mathfrak{sl}_{n}$ simple roots. Thus the set $\Omega$ with the corresponding weight multiplicities is the weight diagram $\mathcal{N}(\gamma)$ of the irreducible $\mathfrak{sl}_{n}$ representation with the highest weight $\gamma$, shifted by the vector $p$ with the coordinates
\begin{equation}
\label{p}
p_{k}=N-\frac{2s_1+...+2s_n}{n}, \;\;\; k=1, ... , n.
\end{equation}
We can imagine this as a slice of an n-dimensional hypercube, which is normal to $(1, ... , 1)$ and has the center in $p$.

The highest weight $\gamma$ is the projection of $\lambda$ to the dual space $h^{*}_{\mathfrak{sl}_{n}}$ of the Cartan subalgebra of $\mathfrak{sl}_{n}$ and it has the coordinates:
\begin{equation}
\label{gamma}
\gamma_i=\frac{2s_1+...+2s_n}{n}-2s_i.
\end{equation}  
So  we see that the leading asymptotic of $\chi_{\lambda}(e^{t})$  is a character of the $\mathfrak{sl}_n$ irreducible representation with the highest weight $\gamma$, shifted by $p$:
\begin{equation}
  \label{eq:22}
 \chi_{\lambda}(e^t)=\chi^{\mathfrak{sl}_{n}}_\gamma(e^\tau)\cdot e^{(t_1+t_2+...+t_n)\cdot (N-\frac{2s_1+...+2s_n}{n}) } \left(1+\mathcal{O}\left(\frac{1}{N}\right)\right).
\end{equation}
Here $\tau$ is the item in Cartan subalgebra of $\mathfrak{sl}_{n}$ which is the projection of the item  $t\in \mathfrak{h}$ in Cartan subalgebra of $\mathfrak{so}_{2n+1}$. Its coordinates are 
\begin{equation}
  \label{eq:23}
\tau_i=-t_i+t_{i+1},\;\;\; i=1, ... ,n-1
\end{equation}
Changing the parameters from  $\{t_{i}\}$ to $\{\Theta_{i}\}$ we obtain:
\begin{equation}
  \label{eq:24}
\chi_\lambda(e^t)=\chi_\gamma(e^\tau)\cdot e^{(t_1+t_2+...+t_n)\cdot N}\cdot \left(\frac{\Theta_1}{N}\cdot ... \cdot \frac{\Theta_n}{N}\right)^{\frac{s_1+...+s_n}{n}} \cdot\left(1+\mathcal{O}\left(\frac{1}{N}\right)\right).
\end{equation}

\subsection{Final}
We combine the equations \eqref{eq:15},\eqref{eq:18},\eqref{eq:24}  together and obtain:
\begin{multline}
  \label{eq:21}
P_\lambda^N(t)=\frac{M_\lambda^{\omega_n,N} \chi_\lambda(e^t)}{(\chi_\omega(e^t))^N}= \prod_{k=1}^{n} \frac{  N^{s_k}}{(s_{k}+k-1)!} \cdot\prod_{i<j}^n(-s_i-i+s_j+j) \cdot\\
\cdot \left(\frac{\Theta_1}{N}\cdot ... \cdot \frac{\Theta_n}{N}\right)^{\frac{s_1+...+s_n}{n}} e^{-\left(\Theta_1+...+\Theta_n\right)\left(1+\mathcal{O}\left(\frac{1}{N}\right)\right) } \chi_\gamma(e^\tau)\left(1+\mathcal{O}\left(\frac{1}{N}\right)\right).
\end{multline}
Rearranging the factors we arrive to the asymptotic of character measure:
\begin{equation}
  \label{eq:26}
  P_\lambda^N(t)=\prod_{i<j}^n(-s_i-i+s_j+j) \cdot  \prod_{k=1}^{n} \frac{  N^{s_k}  e^{-\Theta_k\left(1+\mathcal{O}\left(\frac{1}{N}\right)\right)}}{(s_{k}+k-1)!} \cdot  \frac{\left(\Theta_1\dots \Theta_n \right)^{\frac{s_1+...+s_n}{n}}}{N^{s_1+s_2+...+s_n}}\chi_\gamma(e^\tau)\left(1+\mathcal{O}\left(\frac{1}{N}\right)\right).
\end{equation}
Taking the limit $N\rightarrow\infty$, we complete the proof of Theorem 1:
\begin{equation}
  \label{eq:25}
P_s(\Theta)=\prod_{i<j}^n(-s_i-i+s_j+j) \cdot  \prod_{k=1}^{n} \frac{  \Theta_k^{\frac{s_1+...+s_n}{n}} \cdot e^{-\Theta_k}}{(s_{k}+k-1)!}\chi_\gamma(e^\tau) .
\end{equation}

\section*{Conclusion}

This result can be interpreted as a random walk with a critical drift  \cite{postnova2020multiplicities}. Connection to random walks is easy to understand in the $\mathfrak{sl}_2$-case for the probability measure defined on the weight diagram of $N$-th tensor power of fundamental representation $L^{\omega}$: $p(\lambda)=\frac{\mathrm{dim}V_{\lambda}}{(\mathrm{dim}L^{\omega})^N}$,
where $\mathrm{dim}V_{\lambda}$ is the dimension of the weight subspace in the representation $(L^{\omega})^{\otimes N}$. In this case the probability distribution $p(\lambda)$ is the binomial distribution and $\mathrm{dim} V_{\lambda}$ is the number of random walks of $N$ steps on weight diagram ending in $\lambda$. Plancherel-type measures can be obtained from such a measure by taking an alternating sign combination of random walks \cite{tate2003lattice}. The character measure has the expected value of the weight $\mathbb{E}[\lambda]=N\eta$, which is then interpreted as the expectation value of a random walk with a drift $\eta$, that is determined by the coordinates $\{t_i\}$. Our result appears in the critical drift regime $t_i\sim \ln N$, where random walks stay near the boundary and thus it generalizes Poisson distribution to the case of random walks on the weight diagram of $\mathfrak{so}_{2n+1}$.

In the upcoming publications we are going to consider a similar asymptotic for the $\mathfrak{sl}_n$-case and search for a general formula for all four series of classical Lie algebras: $\mathfrak{sl}_n, \mathfrak{so}_{2n+1}, \mathfrak{sp}_{2n},\mathfrak{so}_{2n}$.

\section*{Acknowledgements}
This research is supported by RFBR grant No. 18-01-00916. We thank Pavel Nikitin for the useful discussions and suggestions.  

\bibliography{bib}
\end{document}